\documentclass[12pt,a4paper]{article}

\newtheorem{theorem}{Theorem}[section]
\newtheorem{lemma}[theorem]{Lemma}
\newtheorem{proposition}[theorem]{Proposition}

\newtheorem{example}[theorem]{Example}
\newtheorem{remark}[theorem]{Remark}
\newtheorem{corollary}[theorem]{Corollary}

\newtheorem{notations}[theorem]{Notations}
\newcommand{\epf}{\ensuremath{\diamondsuit}}
\newcommand{\ic}{\ensuremath{\mathcal{I}}}
\newcommand{\oc}{\ensuremath{\mathcal{O}}}
\newcommand{\fc}{\ensuremath{\mathcal{F}}}

\newcommand{\ec}{\ensuremath{\mathcal{E}}}
\newcommand{\hc}{\ensuremath{\mathcal{H}}}
\newcommand{\cc}{\ensuremath{\mathcal{C}}}
\newcommand{\lc}{\ensuremath{\mathcal{L}}}
\newcommand{\Sc}{\ensuremath{\mathcal{S}}}
\newcommand{\Pt}{\textbf{P}^3}

\newcommand{\Pq}{\textbf {P}^4}
\newcommand{\Pcq}{\textbf{P}^5}
\newcommand{\Psx}{\textbf{P}^6}
\newcommand{\Pn}{\textbf{P}^n}

\begin{document}

\title{On codimension two subvarieties of $\Psx$.}

\author{Ph. Ellia and D. Franco\thanks{Both authors are partially supported by MURST and Ferrara Univ. in the
framework of the project: "Geometria algebrica, algebra commutativa e aspetti computazionali"}}
\date{September 1, 1999} 
 
\maketitle

\section{Introduction}

One of the most attractive problems in algebraic geometry is Hartshorne's conjecture (\cite{Ha}): "let $X \subset \Pn (\textbf{C})$ be a smooth subvariety, if $dim(X)>{2\over 3}n$ then $X$ is a complete intersection". Due to the connection with the existence of rank two vector bundles, the codimension two case is particularly interesting. Thanks to Barth's result (\cite{Barth}) and since no indecomposable rank two vector bundle on $\Pn$, $n \geq 5$, is known, it is generally believed that any smooth, codimension two subvariety of $\Pn$, $n \geq 6$, is a complete intersection.
\par
In the last twentyfive years there have been some results on this conjecture (e.g. \cite{BVdV}, \cite{Bar}, \cite{Ran}, \cite{Ball_Ch}, \cite{LVdV}, \cite{HS}, \cite{Ho}) which may be summarized as follows: if $e \leq n+1$ or if $d < (n-1)(n+5)$ or if $s\leq n-2$, then $X$ is a complete intersection (here $\omega _X \simeq \oc _X(e)$, $d$ is the degree of $X$ and $s$ is the minimal degree of an hypersurface containing $X$).
\par
These inequalities are more or less direct consequences of the following fact observed by Z. Ran (\cite{Ran}): "Let $X\subset \Pn$ be a codimension two subcanonical subvariety. Set $e(k)=k^2-c_1k+c_2$ where the $c_i$ are the Chern classes of the rank two associated vector bundle. If $k\leq n-2$ and if $e(0)...e(k)\neq 0$, then there exists a $(k+1)$-secant line to $X$ through a general point of $\Pn$; in particular $h^0(\ic _X(k))=0$". It seems difficult to extend this approach further and, indeed, there have been no new result in the last ten years.
\par
Clearly it is enough to prove the conjecture for $n=6$ (or for subcanonical smooth threefolds in $\Pcq$). The main results of this paper are (see Theorem \ref{thm}, Theorem \ref{thm2}, Theorem \ref{P5}):
\vfill
\par
\begin{theorem} Let $X \subset \Psx$ be a smooth, codimension two subvariety, if $s\leq 5$ or if $d \leq 73$, then $X$ is a complete intersection.
\par
Let $X \subset \Pcq$ be a smooth, subcanonical threefold. If $s \leq 4$, then $X$ is a complete intersection.
\end{theorem}

In the case $5\leq n\leq 6$ this improves \cite{Ran}, \cite{HS}, \cite{Ho}. We also prove some partial results on "numerical complete intersections" (see  Cor. \ref{nci}, Prop. \ref{Pnci}).
\par
The starting point of our investigation is the following remark: Lefschetz's theorem asserts
that if $X$ is a Cartier divisor on the hypersurface $\Sigma\subset \Pn$, $n\geq 4$, then $X$ is the complete intersection 
of $\Sigma$ with another hypersurface; hence, if $X$ is not a complete intersection $\Sigma$ has
to be singular, more precisely we must have $X \cap Sing(\Sigma) \neq \emptyset$. In fact
something more precise is true: if $X$ is not a complete intersection, then $dim(X \cap Sing(\Sigma))
\geq n-4$: indeed otherwise, cutting with a general $\Pq$ we would get a smooth surface $S$ 
lying on an hypersurface $\Sc$ and with $S \cap Sing(\Sc) = \emptyset$; by Severi-Lefschetz theorem
it would follow that $S$ is a complete intersection, and this in turn would imply that $X$ itself is a 
complete intersection. So we focus our attention on the singular locus of $\Sigma$, an hypersurface of minimal degree containing $X$. More precisely, when $n=6$, $s=5$ and since $e>0$, $X$ is bilinked to a non-reduced subscheme $Z$ on the hyperquintic $\Sigma$ and we show, under suitable assumptions, that $Z_{red} \subset Sing(\Sigma)$ (this last fact follows from a result, see Prop. \ref{pmu}, Prop. \ref{pmuc}, on multiple structures on space curves which might be of independent interest). Combining this with various other considerations and making extensive use of Zak's results (\cite{Zak}, \cite{Zak2}), we conclude if $s=5$. The same approach works also for $n=5$, $s=4$ and in some cases, if $n=6$, $s>5$ (see Theorem \ref{P5}, Theorem \ref{thm2}).
\par
We thank Steven Kleiman and Christian Peskine for useful comments.   

\section{Generalities.}

\emph{Notations.}
\par
In this section, $X \subset \Pn, n \geq 4$, will denote a smooth, non degenerated, codimension two
subvariety of degree $d$. We will often (but not always) assume $X$ subcanonical: $\omega _X \simeq
\oc _X(e)$; notice that, thanks to Barth's theorem, this condition is automatically fullfilled if
$n \geq 6$.
\par
As usual we will denote by $s(X)$ (or more simply $s$) the minimal degree of an hypersurface
containing $X$.
\bigskip

As explained in the introduction, we have:

\begin{lemma}
\label{l1}
Let $X \subset \Pn$, $n \geq 4$, be a smooth codimension two subvariety. Assume $X$
lies on the hypersurface $\Sigma$. If $X$ is not a complete intersection, then $dim(X \cap
Sing(\Sigma)) \geq n-4$.
\par \noindent
Moreover if $deg(\Sigma )=s$, then $n-4 \leq dim(Sing(\Sigma) \cap X) \leq n-3$.
\end{lemma}

For the last assertion of the lemma notice that, by minimality of $s$, $X \not\subset Sing(\Sigma)$.
\par 
In the sequel we will concentrate on the case where $\Sigma$ has minimal degree $s$ and we will
show that under suitable assumptions $dim(X \cap Sing(\Sigma))=n-4$. The first part of the next
lemma is just a reformulation of a result of Ellingsrud-Peskine (\cite{EP}):

\begin{lemma}
\label{lsd} 
\begin{enumerate}
\item Let $X \subset \Pn , n \geq 4$, be a smooth codimension two subvariety of degree $d$ 
with $\omega _X \simeq \oc _X(e)$, then: 
\par \noindent 
$s(n+1+e)-s^2 \leq d \leq s(n-1+e)+1$.
\item Assume moreover $Pic(X) \simeq \textbf {Z}.H$. Let $\Sigma$ be an hypersurface of degree $s$ 
containing $X$ and suppose $dim(Sing (\Sigma)\cap X) = n-3$. Let $T \sim lH$ denote the divisorial
part of $Sing(\Sigma) \cap X$, then:
\par \noindent
a) $(s-l)(n+1+e)-(s-l)^2 \leq d \leq (s-l)(n-1+e)+1$, 
\par \noindent
b) $h^1(\ic _X(l))\neq 0$.
\end{enumerate}
\end{lemma}

\textbf{Proof.} 
\begin{enumerate}
\item Apply Lemme 1 of \cite{EP} to a section of $X$ with a general $\Pq$. 
\item
(a) This is just a slight variation on the first point. For the convenience of the reader we will
sketch it briefly. The inclusion $X \subset \Sigma$ induces $\sigma :\oc _X \to N_X^*(s)$. The zero 
locus $(\sigma )_0$ is the scheme theoretical intersection of the jacobian of $\Sigma$ with $X$. 
By hypotesis, $(\sigma )_0$ has codimension one in $X$ and we can divide $(\sigma )_0$ 
by the codimension one part getting a section $\sigma '\in H^0(N_X ^*(s-l))$ thus $(\sigma ')_0 = 
c_2(N_X^*(s-l))$ in $A^2(X)$. A short calculation, using $c_1(N_X) = (n+1)H+K = (n+1+e)H$ and 
$c_2(N_X) = dH^2$ ("formule clef"), gives: 
$c_2(N_X^*(s-l))= c_2(N_X (-s+l)) = [-(s-l)(n+1+e)+d+(s-l)^2].H^2$.
Taking degree: 
$deg(\sigma ')_0 = [-(s-l)(n+1+e)+d+(s-l)^2].H^{n-2} = 
[-(s-l)(n+1+e)+d+(s-l)^2]d$. It follows that 
$d \geq (s-l)(n+e+1)-(s-l)^2$. 
On the other hand, if $F=0$ is an equation of $\Sigma$, the jacobian is contained in 
${\bf V}(F'_i) \cap {\bf V}(F'_k) \cap X$, for suitable indices, and 
the linear system cut out by the partials of $F$ has a degree $l$
fixed part and a degree $s-l-1$ moving one, hence $(\sigma ')_0$
is contained in the intersection of two divisors
of degree $s-l-1$ without common components. 
This implies that $(s-l-1)^2.H^2 - [-(s-l)(n+1+e)+d+(s-l)^2].H^2$ is effective. 
Taking degrees again yields: $d \leq (s-l)(n-1+e)+1.$
\par
(b) Suppose by contradiction that $X$ is $l$-normal.
Since $Pic(X) \simeq {\bf Z}.H$, $T$ is cut out on $X$ by an hypersurface: $T = X \cap {\bf V} (P)$, 
where $P$ is an homogeneous polynomial. Let $F = 0$ be an equation of $\Sigma$ and set $F'_i = 
\partial F /\partial x_i$. Since $T \subset (Sing(\Sigma)\cap X)$, we have $X\cap {\bf V}(P) 
\subset X \cap {\bf V}(F'_i)$. This implies that $P$ divides $F'_i, mod({\bf I}(X)$. So: 
$PG_i -F'_i \in {\bf I}(X)$. By minimality of $s$, this implies $PG_i = F'_i$. By Euler's relation 
($ch(k) = 0$), $P$ divides $F$, but this is impossible, since by minimality of $s$ again, 
$\Sigma$ is integral. \epf
\end{enumerate}
\bigskip

By the way, the previous lemma has an interesting consequence:

\begin{corollary}
\label{nci}
Let $X \subset \Pn$, $n \geq 4$, be a smooth, codimension two subvariety. 
Assume $X$ is numerically a complete intersection $(a,b), a \leq b$ (i.e. $d = ab, \omega _X \simeq 
\oc _X(a+b-n-1)$).
\begin{enumerate}
\item If $b > a(a-3) + 3$ then $X$ is a complete intersection.
\item If $a \leq n-1$ then $X$ is a complete intersection.
\end{enumerate}
\end{corollary}

\textbf{Proof.}
(a) Let $C \subset \Pt$ be a general space section of $X$. Then $C$ has the numerical characters of a 
complete intersection $(a,b)$ in $\Pt$, in particular $\pi = g(C) = G(d,a)$ (recall that $G(d,s)$ 
is the maximal genus of a smooth curve of degree $d$ not contained in a surface of degree $<s$). 
If $s(C) \geq a$ then, by \cite{GP}, $C$ is a complete intersection. It follows that $X$ also is a 
complete intersection. So we may assume $s(C) < a$. By Lemma \ref{lsd}: $d \leq s(n-1+e)+1$, this 
implies: $ab \leq (a-1)(a+b-2)+1$ and the result follows.
\par \noindent
(b) Arguing as above if $s(X) \geq n-1$ we are done. If $s(X) < n-1$ we conclude with \cite{Ran}. \epf

\begin{remark} 
\begin{enumerate}
\item We recall that Ran's theorem \cite{Ran} implies that a numerical complete intersection 
$(a,b)$ with $a \leq n-2$ is a complete intersection.
\item The second statement of the corollary can be found in \cite{HS}. 
\item Another immediate consequence of Lemma \ref{lsd} is: $e \geq -n+1$ (cp with \cite{Ball_Ch})
\end{enumerate}
\end{remark}

\section{Codimension two subvarieties lying on quintic hypersurfaces in $\Psx$.}

\subsection{Preliminaries} 
\begin{notations}
\label{not2}
In this section $X \subset \Psx$ will denote a smooth, codimension two 
subvariety, of degree
$d$, with $\omega _X \simeq \oc _X(e)$ and with $s(X)=5$. We will assume that $X$ is not a complete
intersection and derive a contradiction.
\par
We will denote by $\Sigma$ an 
irreducible quintic hypersurface containing $X$. (Notice that we may assume $d>25$ and hence
$\Sigma$ uniquely determined.)
\par 
By Serre's construction we may associate to $X$ a rank two vector bundle:
$$ 0 \to \oc \to E \to \ic _X(e+7) \to 0 $$
The Chern classes of $E$ are: $c_1(E)=e+7,c_2(E)=d$.
\par
Since $h^0(\ic _X(5))=1$, $E(-e-2)$ has a section and this is the least twist of $E$ having a
section. If $X$ is not a complete intersection, $E$ doesn't split and the section of $E(-e-2)$
vanishes in codimension two:
$$ 0 \to \oc \to E(-e-2) \to \ic _Z(-e+3) \to 0$$
where $Z \subset \Psx$ is a locally complete intersection subscheme of degree $d(Z)=c_2(E(-e-2))=
d-5e-10$ and with $\omega _Z \simeq \oc _Z(-e-4)$.
\par
Let $Y$ be a section of $Z$ with a general $\Pt$, then $Y \subset \Pt$ is a locally complete 
intersection curve of degree $d(Z)$, with $\omega _Y \simeq \oc _Y(-e-1)$, in particular
$2p_a(Y)-2=-d(Z)(e+1)$. Let $Y_{red}=Y_1 \cup ... \cup Y_r$ be the decomposition into irreducible
components. Making a primary decomposition we can write $Y = \tilde{Y}_1 \cup ... \cup \tilde{Y}_r$
where $\tilde{Y}_i$ is a locally Cohen-Macaulay (and generically l.c.i.) subscheme with support
$Y_i$. So $\tilde{Y}_i$ is a locally Cohen-Macaulay multiple structure of multiplicity $m_i$
on the integral curve $Y_i$. The multiplicity $m_i$ is determined by: $deg(\tilde{Y}_i)=m_i.deg(Y_i)$.
\end{notations}

\begin{lemma}
\label{mi>1}
With notations as above $Z$ is non-reduced, more precisely $m_i \geq 2$, $\forall i,1\leq i \leq r$.
\end{lemma}

\textbf{Proof.}
First we claim that $\omega _{\tilde{Y}_i}\simeq \omega _{Y}\otimes \ic _{\Delta ,\tilde{Y}_i}$,
where $\Delta $ is the scheme theoretic intersection of $\tilde{Y}_i$ and
$W_i:= \cup _{j\not = i} \tilde{Y}_j$
\par\noindent
Indeed, by the following sequence:
$$
0\to \ic _Y \to \ic _{W_i} \to \omega _{\tilde{Y}_i}\otimes \omega _{Y}^{-1} \to 0
$$
we get 
$$
\omega _{\tilde{Y}_i}\simeq {\ic _{W_i}\over \ic _Y}\otimes \omega _Y
\simeq {\ic _{W_i}+\ic _{\tilde{Y}_i}\over \ic _{\tilde{Y}_i}}\otimes \omega _{Y}\simeq
\omega _{Y}\otimes \ic _{\Delta ,\tilde{Y}_i}
$$ 
and the claim follows.
\par
If $m_1 = 1$ then $\omega_Y|Y_1 \simeq \oc _{Y_1} (-e-1)$, but we also have 
$\omega_{Y_1} \simeq \omega _{Y}\otimes \ic _{\Delta ,Y_1} \simeq \ic _{\Delta,Y_1}(-e-1)$, hence $\chi (\omega _{Y_1})\leq \chi (\oc _{Y_1}(-e-1))$.
Since the arithmetic genus of an integral curve is positive, we get a contradiction. \epf
\bigskip

The next lemma controls the dimension of $X \cap Sing(\Sigma)$ (and also of $X\cap Z$).

\begin{lemma}
\label{l2}
With notations as \ref{not2} we have:
\begin{enumerate}
\item $dim(Sing(\Sigma) \cap X)=2$.
\item $Z$ is contained in $\Sigma$ and $dim(X \cap Z)=2$.
\end{enumerate}
\end{lemma}

\textbf{Proof.} 
(a) By Lemma \ref{l1}, $2 \leq dim(Sing(\Sigma) \cap X) \leq 3$. Suppose $dim(Sing(\Sigma) \cap X)
=3$ and let $T \sim lH$ denote the divisorial part of $Sing(\Sigma) \cap X$. By Lemma \ref{lsd},
$d \leq (5-l)(e+5)+1$. By Zak's theorem (\cite{Zak}), $h^1(\ic _X(1))=0$, and, again by Lemma \ref{lsd}, we
may assume $l \geq 2$. It follows that $d \leq 3e+16$. This implies (see Notations \ref{not2}):
$d(Z) \leq -2e+6$, since we may assume $e \geq 8$ (see \cite{HS}), we get a contradiction.
\par \noindent
(b) By construction $X$ and $Z$ are bilinked on $\Sigma$, hence $Z \subset \Sigma$. The inclusion
$X \subset \Sigma$ induces a section $\oc _X \to N_X^*(5)$ whose zero locus is $Jac(\Sigma) \cap
X$. Since $N_X^*(5) = E(-e-2)_{|X}$, this section is nothing else than the restriction to $X$ of the
section of $E(-e-2)$ vanishing along $Z$. It follows that $X \cap Z=Jac(\Sigma) \cap X$ schematically.
By (a) we conclude that $dim(X \cap Z)=2$. \epf
\bigskip

We will go on, step by step, proving that:
\begin{itemize}
\item $Z_{red}$ doesn't contain any irreducible component of degree one,
\item $Z_{red}$ doesn't contain any irreducible component of degree two,
\item $Z_{red}$ doesn't contain any irreducible component of degree greater or equal to three.
\end{itemize}

Then it will follow that $X$ has to be a complete intersection.

\subsection{$Z_{red}$ doesn't contain any irreducible component of degree one.}
\bigskip

\begin{proposition}
\label{compd=1}
With notations as in \ref{not2}, $Z_{red}$ doesn't contain any irreducible component of degree one.
\end{proposition}

\textbf{Proof.}
Let $L \subset Z_{red}$ be a codimension two linear subspace. Consider the linear system, $\delta$,
cut out on $X$ by the hyperplanes through $L$. By Bertini's theorem, the general member, $V$, of
$\delta$ is smooth outside the base locus $B_{\delta}=L \cap X$. Now let $x \in L\cap X$ such 
that every divisor of $\delta$ is singular at $x$. If $V$ is such a divisor, then $dim(T_xV) \geq
4$. Since $T_xV=T_xX \cap H$ it follows that $T_xX \subset H$. If this is true for every hyperplane
$H$ containing $L$, we have $T_xX=L$. Now consider the Gauss map $g:X \to Gr(4,6)$ defined by
$g(x)=T_xX$. By (another) theorem of Zak (\cite{Zak2}), $g$ is finite. It follows that if $V$ is 
sufficiently
general in $\delta$, then $dim(Sing(V))\leq 0$. In particular $V \subset H \simeq \Pcq$ is irreducible
(two threefolds in $\Pcq$ intersect at least along a curve). Now $V=X \cap H \subset \Sigma \cap H=
L \cup F$ where $F$ is a quartic hypersurface of $\Pcq \simeq H$. Since $V$ is irreducible, we have
$V \subset L$ or $V \subset F$. 
\par
If $V \subset L$ then we are done. (Consider the exact sequence:
$ 0 \to \ic _X \to \ic _X(1) \to \ic _V(1) \to 0$.)
\par
So we may assume $s(V) \leq 4$. Let $S \subset \Pq$ be a general hyperplane section of $V$. Then
$S$ is a smooth (because $dim(Sing(V)) \leq 0$), degree $d$, surface in $\Pq$ with $s(S) \leq 4$. 
By Lemma \ref{lsd}: $d \leq 4(e+5)
+1$ (observe that the quantity $e+n$ is invariant by hyperplane section). It follows that
$d(Z) = d-5e-10 \leq -e+11$. Since we may assume $e \geq 8$ (\cite{HS}), we get $d(Z) \leq 3$.
\par
By Lemma \ref{mi>1}, the only possibility is that $Z$ is a locally complete intersection multiple
structure on $L$ of multiplicity $r$, $2 \leq r  \leq 3$. By \cite{Ma}, any such multiple structure
in $\Psx$ is a complete intersection, so $Z$ is a complete intersection. This implies that $E$ splits
and that $X$ too is a complete intersection, contradiction. \epf

\subsection{$Z_{red}$ doesn't contain any irreducible component of degree two.}

We begin with a lemma which will be useful also in other circumstances:

\begin{lemma}
\label{Clin}
Let $X \subset \Psx$ be a smooth, codimension two subvariety. Let $V=X\cap H$ be an hyperplane
section of $X$. If $C \subset \Pt$ is a section of $V$ with a general $\Pt$, then $C$ is a
smooth, irreducible curve which is linearly normal in $\Pt$.
\end{lemma}

\textbf{Proof.}
The threefold $V$ will be singular at the points where the hyperplane $H$ is tangent to $X$.
By Zak's theorem on tangencies (see \cite{LVdV}, p.18), $dim(Sing(V)) \leq 1$. It follows that
the intersection of $V$ with a general $\Pt$ is a smooth curve $C \subset \Pt$.
\par
Now we are going to show that $C$ is linearly normal in $\Pt$ (i.e. $h^1(\ic _C(1))=0$). Consider
the following exact sequences of restriction to an hyperplane:
\begin{equation}
 0 \to \ic _X(m-1) \to \ic _X(m) \to \ic _V(m) \to 0    \label{eq:Xm}
\end{equation}
\begin{equation}
 0 \to \ic _V(m-1) \to \ic _V(m) \to \ic _S(m) \to 0    \label{eq:Vm}
\end{equation}
\begin{equation}
 0 \to \ic _S(m-1) \to \ic _S(m) \to \ic _C(m) \to 0    \label{eq:Sm}
\end{equation}
By (\ref{eq:Sm}) for $m=1$, we see that $h^1(\ic _S(1))=h^2(\ic _S)=0$ implies $h^1(\ic _C(1))=0$.
By (\ref{eq:Vm}) for $m =0,1$, we see that these vanishings follow from
 $h^j(\ic _V(2-j))=0$,
$1 \leq j \leq 3$. Finally, by (\ref{eq:Xm}) for $m=-1,0,1$, we see that $h^j(\ic _X(2-j))=0$,
$1 \leq j \leq 4$ gives the result. We have $h^j(\ic _X(2-j))=h^{j-1}(\oc _X(2-j))=0$, if $3 \leq
 j \leq 4$, by Kodaira's theorem. By Barth's theorem $h^2(\ic _X)=h^1(\oc _X)=0$; finally by
Zak's theorem $h^1(\ic _X(1))=0$ and we are done.
\par
Since $C$ is smooth, $C$ is irreducible if and only if it is connected. If $C$ is not connected
it is the disjoint union of several smooth curves: $C = C_1\cup ...\cup C_t$. Since $C$ is
linearly normal, $h^0(\oc _C(1))=4$, but $h^0(\oc _C(1))= \Sigma h^0(\oc _{C_i}(1))$, and we see
that the only possibility is $t=2$ and $C_i$ is a line, $1 \leq i \leq 2$. This implies that
$S \subset \Pq$ is the union of two planes meeting at one point, but this is impossible since
such a surface is not locally Cohen-Macaulay (\cite{Ha2}) (and a fortiori not l.c.i.). \epf

\begin{proposition}
\label{compd=2}
With notations as in \ref{not2}, $Z_{red}$ doesn't contain any irreducible component of
degree two.
\end{proposition}

\textbf{Proof.}
Suppose $Q$ is an irreducible component of $Z_{red}$ of degree two. Of course, $Q$ is degenerated
in $\Psx$, denote by $H$ the hyperplane containing $Q$. Now by Zak's theorem on tangencies,
$V=X \cap H$ has a singular locus of dimension at most one. We have $V=X \cap H \subset \Sigma \cap
H=Q \cup F$ where $F$ is a cubic hypersurface of $H \simeq \Pcq$. Since every irreducible component 
of $V$ has dimension three and since $dim(X \cap Z)=2$ (see Lemma \ref{l2}), we get $V \subset F$. 
Now if $S$ is a general hyperplane section of $V$, we have $h^0(\ic _S(3)) \neq 0$. If 
$h^0(\ic _S(2))\neq 0$ then, $C$, a general hyperplane section of $S$ lies on a surface of degree
$\leq 2$, since $C$ is subcanonical, smooth and irreducible (see Lemma \ref{Clin}) this implies that
$C$ is a complete intersection. It follows that $X$ is a complete intersection too, contradiction.
Hence we may assume that $S$ lies on an irreducible cubic hypersurface, $F_H$.
\par
We claim that $F_H$ is a normal cubic hypersurface of $\Pq$. Indeed otherwise, $C$, a general
hyperplane section of $S$ would lie on a cubic surface with a double line. Such a cubic
surface is the projection of a cubic scroll in $\Pq$ and we would have $h^0(\oc _C(1)) > 4$ 
contradicting Lemma \ref{Clin}.
\par
Now we conclude with \cite{Koe}, Thm.4.1 that $S$ is a complete intersection and this
yields the desired contradiction. \epf

\subsection{$Z_{red}$ doesn't contain any irreducible component of degree $\geq 3$.}

Thanks to the previous results we may assume that every irreducible component of $Z_{red}$ has 
degree at least three. Our first task will be to show that under this condition $Z_{red}$ is contained in $Sing(\Sigma)$ (notations are as in \ref{not2}). The proof of this fact will follow from a
general result about multiple structures on space curves.
\par
Once we will have proved that $Z_{red}$ is contained in $Sing(\Sigma)$ we will conclude the proof
distinguishing two subcases: I) $Z_{red}$ contains an irreducible component of degree three; 
II) every irreducible component of $Z_{red}$ has degree at least four.

\subsubsection{Multiple structures on space curves.}

Let $C \subset S \subset \Pt$ be an integral curve lying on a smooth surface. There exists a
uniquely determined loc. C.M. multiple structure of multiplicity $m$ on $C$, $C_m$, which lies on $S$.
Indeed it is the (Weil hence Cartier) divisor $mC$ on $S$. By adjunction formula: $C(C+K)=
2g-2$, where $g:=p_a(C)$. Since $K \sim (s-4)H$ where $s=deg(S)$, we get $C^2=2g-2-d(s-4)$ (here
$d=deg(C)$).
Now we have $mC(mC+K)=2p-2$ where $p:=p_a(C_m)$ and we easily compute: $p=\mu (d,g,s,m)$ where
\begin{equation}
\mu (d,g,s,m):=1+m^2(g-1)-(s-4)d{m(m-1)\over 2}			\label{eq:mu}
\end{equation}

Observe that the arithmetic genus of $C_m$ doesn't depend on the curve $C$ nor on the surface $S$ but
just on their numerical invariants.
\par
Now assume that $S$ is singular with $dim(C \cap Sing(S))=0$. There is still a uniquely determined
loc.C.M. multiple structure of multiplicity $m$ on $C$, $C_m$, contained in $S$. This structure can
be defined as follows: $C_m$ is the greatest loc.C.M. subscheme of $C^{(m)} \cap S$ where
$C^{(m)}$ denotes the $m$-th infinitesimal neighbourhood of $C$. Now a natural question is: what
can be said on the arithmetic genus of $C_m$? A simple example will suggest the answer: a double
line on a smooth quadric has genus $-1$ while a double line on a quadric cone has genus $0$; more
generally it is well known that, if $C$ is smooth, the singularities of $S$ lying on $C$  increase the degree of the subline bundle of $N_C$ defined by $S$. So we may wonder if the inequality $p_a(C_m)\geq \mu (d,g,s,m)$ holds in general. We will show that this is indeed the case if $C$ is Gorenstein. Unfortunately our proof doesn't extend to the general case, but we can prove a similar inequality for double structures and triple subcanonical structures on an integral curve (see Prop. \ref{pmuc}) and this will be enough for our purposes.

\begin{proposition}
\label{pmu}
Let $C \subset S \subset \Pt$ be an integral, Gorenstein curve of degree $d$, arithmetic genus $g$, lying
on the irreducible surface $S$ of degree $s$. Assume $dim(C \cap Sing(S))\leq 0$ and let
$C_m$ be the unique loc.C.M. multiplicity $m$ structure on $C$ contained in $S$. Then $p_a(C_m)
\geq \mu (d,g,s,m)$.
\end{proposition}

To enlight the proof, observe that the expression of $\mu (d,g,s,m)$ if $S$ is smooth can be found also in the following way: set $\fc _i=\ic _{C_i,C_{i+1}}$. We have $\fc _i \simeq \ic _{C_i,S}\otimes \oc _C \simeq (\omega ^{\surd}_C(s-4))^{\otimes i}$. From the sequence: $0 \to \fc _i \to \oc _{C_{i+1}} \to \oc _{C_i} \to 0$, we get $p_{i+1}=p_i+g-1-deg(\fc _i)$. Now $deg(\fc _i)=deg(\fc _{i-1})+2-2g+d(s-4)$ and by induction, starting from $\fc _0 \simeq \oc _C$, we get the result. If $dim(C\cap Sing(S))=0$ we repeat this argument showing $deg(\fc _i)\leq deg(\fc _{i-1})+2-2g+d(s-4)$, taking into account that $\fc _i \simeq (\ic _{C_i,S}\otimes \oc _C)^{\circ}$ ($(\fc)^{\circ}$ denotes $\fc$ mod. torsion). Recall that a rank one torsion-free sheaf on an integral curve, $C$, is of the form $\fc \simeq \ic _Z(D)$ where $Z$ is a zero-dimensional subscheme of $C$ and where $D$ is a Cartier divisor, the degree of such a sheaf is defined by $-deg(Z)+deg(D)$. From Riemann-Roch for Cartier divisors (\cite{Fu}) it follows that $\chi (\fc ) =deg(\fc )-g+1$. Moreover if $C$ is Gorenstein then every rank one torsion-free sheaf is reflexive and we have $deg(\fc ^{\surd})=-deg(\fc )$ (see \cite{Ha3}); this equality is not always true for reflexive sheaves on a non-Gorenstein curve (see Example \ref{nonGo}) and this is the main obstruction to have a general statement.
\bigskip

\textbf{Proof.}
Applying $\hc om_{\oc _S}(-,\omega _S)$ to the sequence:
$$
0 \to \ic _{C_i,S} \to \oc _S \to \oc _{C_i} \to 0
$$
we get: $0 \to \omega _S \to \hc om_{\oc _S}(\ic _{C_i,S},\omega _S) \to \omega _{C_i} \to 0$. Set $\lc :=\hc om_{\oc _S}(\ic _{C_i,S},\oc _S)$, so that the above sequence reads like: 
$$
0 \to \omega _S \to \lc (s-4) \to \omega _{C_i} \to 0
$$
Restricting to $C$ yields: $\lc \mid _C(s-4) \to \omega_{C_i}\mid _C \to 0$ (*).
\par
Applying $\hc om_{\oc _{C_i}}(-,\omega _{C_i})$ to $0 \to \fc _{i-1} \to \oc _{C_i} \to \oc _{C_{i-1}} \to 0$, gives:
$$
0 \to \omega _{C_{i-1}} \to \omega _{C_i} \to \hc om_{\oc _{C_i}}(\fc _{i-1},\omega _{C_i}) \to 0
$$
Restricting to $C$: $\omega _{C_i}\mid _C \to \hc om_{\oc _C}(\fc _{i-1}, \omega _C) \to 0$ (because $\fc _{i-1}$ is an $\oc _C$-module). Combining with (*) we get a morphism: $\lc \mid _C(s-4) \to \hc om_{\oc _C}(\fc _{i-1},\omega _C) \to 0$. Since $\hc om_{\oc _C}(\fc _{i-1},\omega _C)$ is torsion free this yields: $\psi: (\lc \mid _C)^{\circ}(s-4) \to \hc om_{\oc _C}(\fc _{i-1},\omega _C) \to 0$ ($\circ$). It follows that:
\par
$deg(\hc om_{\oc _C}(\fc _{i-1},\omega _C)) \leq deg((\lc \mid _C)^{\circ}(s-4))$ (**).
\par
On the other hand since $\lc = \hc om_{\oc _S}(\ic _{C_i,S}, \oc _S)$, we have a morphism: $\lc \mid _C \to \hc om_{\oc _C}(\ic _{C_i,S}\mid _C,\oc _C)$. Since $\hc om_{\oc _C}(\ic _{C_i,S}\mid _C,\oc _C)$ is torsion free, we get: $\varphi :(\lc \mid _C)^{\circ}(s-4) \to \hc om_{\oc _C}(\ic _{C_i,S}\mid _C,\oc _C)(s-4)(\circ \circ)$, and this map is clearly injective. Since $\fc _i=(\ic _{C_i,S}\otimes \oc _C)^{\circ}$, $\hc om_{\oc _C}(\ic _{C_i,S}\mid _C,\oc _C)=\hc om_{\oc _C}(\fc _i,\oc _C)=:\fc _i^{\surd}$. In conclusion: $deg(\fc _i^{\surd}(s-4))\geq deg((\lc \mid _C)^{\circ}(s-4)$. Combining with (**): $deg(\hc om_{\oc _C}(\fc _{i-1},\omega _C)) \leq deg(\fc _i^{\surd}(s-4))$. Now since $C$ is Gorenstein, $deg(\fc _i^{\surd})=-deg(\fc _i)$. Since $deg(\hc om_{\oc _C}(\fc _{i-1},\omega _C))=-deg(\fc _{i-1})+2g-2$, we get the desired relation:
\begin{equation}
\label{eq:deg}
deg(\fc _i) \leq deg(\fc _{i-1})-2g+2+d(s-4)
\end{equation}
From $0 \to \fc _i \to \oc _{C_{i+1}} \to \oc _{C_i} \to 0$, we have $p_a(C_{i+1})=p_a(C_i)-\chi (\fc _i)$, by Riemann-Roch, $\chi(\fc _i)=deg(\fc _i)-g+1$ and we conclude. \epf

\begin{remark}
The previous proof has been inspired by \cite{Flo}.
\end{remark}
\bigskip

Let $C$ be an integral curve and denote by $p: \tilde C \to C$ the normalization. As usual we put: $\delta = h^0(p_*\oc _{\tilde C}/\oc _C)$.

\begin{proposition}
\label{pmuc}
With notations as in Proposition \ref{pmu} (but $C$ non necessarily Gorenstein), assume one of the following holds: $m=2$ or $m=3$ and $C_3$ subcanonical, then $p_a(C_m) \geq \mu (d,g,s,m)-{m(m-1)(2\delta -1)\over 2}$.
\end{proposition}
\bigskip

The proof rests on the following:
\begin{lemma}
\label{omd}
Let $C$ be an integral curve of arithmetic genus $g$, then:
\par
$deg(\hc om_{\oc _C}(\omega _C, \oc _C)) \leq -2g+2+2\delta -1$.
\end{lemma}

\textbf{Proof.}
By the projection formula $p_*(\omega _{\tilde C}) \simeq p_*(\oc _{\tilde C}) \otimes L$, here $L$ is invertible of degree $2g'-2$, $g'=g-\delta$. Dualizing the injection $p_*\omega _{\tilde C} \hookrightarrow \omega _C$, we get:
$$
\hc om_{\oc _C}(\omega _C,\oc _C) \hookrightarrow \cc \otimes L^{\surd}
$$
where $\cc =\hc om_{\oc _C}(p_* \oc _{\tilde C},\oc _C)$ is the conductor of $p$.
\par 
Hence $deg(\hc om_{\oc _C}(\omega _C,\oc _C)) \leq deg(\cc)-2g'+2$. Recall that $\cc$ is also an $\oc _{\tilde C}$-ideal sheaf. Denote by $n$ the degree of the subscheme of $\tilde C$ defined by $\cc$. We claim that: $\delta + \gamma =n$ where $\gamma = deg(\Gamma)$. Indeed this is a local question, so let $A$ be a one-dimensional integral local ring, $A'$ its integral closure and $I$ the conductor; then the claim follows from ${(A'/I)\over (A/I)} \simeq {A'\over A}$. Since $deg(\cc)=-\gamma =-n + \delta$ and taking into account that $\delta +1 \leq n$ (\cite{Serre}, p.80), we get the result. \epf

\pagebreak

\textbf{Proof of Prop. \ref{pmuc}.}
\par
Dualizing the morphisms $\psi$, $\varphi$ (see $(\circ)$ and $(\circ \circ)$ in the proof of Prop. \ref{pmu}) and using that $\fc _1$ injects in its bidual, we get $\fc _1 \hookrightarrow \hc om_{\oc _C}(\omega _C,\oc _C)(s-4)$; the case $m=2$ follows from Lemma \ref{omd}.
\par
Again as in Prop. \ref{pmu} we have:
\par
$deg(\hc om_{\oc _C}(\fc _2,\oc _C))(s-4) \geq deg(\hc om_{\oc _C}(\fc _1,\omega _C))$
\par \noindent
$=-deg(\fc _1)+2g-2 \geq 4g-4-d(s-4)-2\delta +1$ (*), where the last follows from the first step. On the other hand, since $\omega _{C_3} \simeq \oc _{C_3}(a)$, we have $\fc _2 \simeq \omega _C(-a)$ hence $deg(\hc om_{\oc _C}(\fc _2,\oc _C)) \leq da -2g+2+2\delta -1$, combining with (*) and after a short computation, the result follows. \epf

\begin{example}
\label{nonGo}
We show that it is not always true that $deg(\fc^{\surd})=-deg\fc$ for $\fc$ a rank one reflexive sheaf on a non-Gorenstein integral curve.
\par
The linear system $\mid C_0+2f \mid$ on $\textbf{F}_2=\textbf{P}(\oc \oplus \oc (-2))$ maps $\textbf{F}_2$ to a quadric cone $g(\textbf{F}_2)=Q\subset \Pt$. Consider a smooth curve $\tilde C \in \mid C_0+5f\mid$, then $C =g(\tilde C)$ is a curve of degree $5$, arithmetic genus $2$ and $g:\tilde C \to C$ is the normalization (so $\delta = 2$). The curve $C$ is not Gorenstein: $C$ is a linked to a line $L$ by a complete intersection $(Q,F)$, we have $\ic _{L,Q} \simeq \omega _C(-1)$; if $\omega _C$ were invertible , $F$ would have to be a minimal generator of $\ic _{L,Q}$ at the vertex $v$ of the cone $Q$, but this is impossible since $F$ is singular at $v$. With notations as in \cite{Serre}, p. 80, we have $\delta +1 \leq n \leq 2\delta$ and $n \neq 2\delta$ since $C$ is not Gorenstein. It follows that $n=3$ and $\gamma =1$. We see that the conductor $\cc$ defines the point $v$ on $C$. So $deg(\cc)=-1$. Observe that $\cc = \hc om_C(g_*\oc _{\tilde C}, \oc _C)$ is reflexive. Now from the exact sequence $0 \to \oc _C \to \cc^{\surd} \to \ec xt^1_{\oc _C}(\oc _{C,v},\oc _C) \to 0$ we get $deg(\cc ^{\surd})\neq 1$. Indeed since $\oc _{C,v}=A$ is not Gorenstein, $dim_k Ext^1_A(k,A)\neq 1$ (in fact it is equal to $2$). 
\end{example}

\subsubsection{$Z_{red}$ is contained in $Sing(\Sigma)$.}

We recall the notations of \ref{not2}: if $Y \subset \Pt$ denotes the intersection of $Z$ with
a general $\Pt$ then $\omega _Y \simeq \oc _Y(-e-1)$ and we can write $Y= \tilde{Y}_1 \cup ...\cup 
\tilde{Y}_r$ where $\tilde{Y}_i$ is a multiplicity $m_i$-structure on the integral curve $Y_i$.
We set $g_i=p_a(Y_i)$ and $d_i=deg(Y_i)$ so that $d(Z)=deg(Y)=\Sigma m_id_i$.

\begin{proposition}
\label{ZrinSing}
With notations as in \ref{not2}, if every irreducible component of $Z_{red}$ has degree at least
three, then $Z_{red} \subset Sing(\Sigma)$.
\end{proposition}

\textbf{Proof.}
We have $Y \subset F \subset \Pt$ where $F=\Sigma \cap \Pt$ is an irreducible quintic surface.
Clearly it is enough to show that $Y_{red}= Y_1 \cup ...\cup Y_r$ is contained in $Sing(F)$.
\par
First let's observe that $deg(Y)=d(Z)\leq 16$. Indeed by Lemma \ref{lsd}: $d \leq 5e+26$ so since $d(Z)=d-5e-10$ (see \ref{not2}), we get $d(Z) \leq 16$.
\par
Now assume that $deg(Y_i)=d_i\leq 4,\forall i$. Since every integral curve of degree $\leq 4$ is Gorenstein, we can apply Proposition \ref{pmu}:
if $Y_1$ is not contained in $Sing(F)$, we must have: 
$p_1:=p_a(\tilde{Y}_1) \geq \mu (d_1,g_1,5,m_1)$ (*). Since $\omega _{\tilde{Y}_1} \simeq
\omega _Y \otimes \ic _{\Delta, \tilde{Y}_1}$ (see the proof of Lemma \ref{mi>1}) where
$\Delta= \tilde{Y}_1 \cap W_1$ and since $\omega _Y \simeq \oc _Y(-e-1)$, we get $2p_1 - 2
\leq -m_1d_1(e+1)$. Combining with (*):
\begin{equation}
-m_1d_1(e+1) \geq 2m_1^2(g_1-1)-m_1(m_1-1)d_1   \label{eq:a}
\end{equation}
This can be written as:
\begin{equation}
m_1d_1 \geq 2m_1(g_1-1)+d_1(e+2)   \label{eq:b}
\end{equation}
Now if $g_1 \geq 1$, since $d_1 \geq 3$ by assumption, $e \geq 8$
(cf \cite{HS}) and $m_1d_1 \leq \Sigma m_id_i = d(Z) \leq 16$,  we get $16 \geq m_1d_1 \geq 30$ which is absurd. So we may assume $g_1=0$.
Since $16 \geq m_1d_1$, we get: $2m_1 \geq d_1(e+2)-16$. Since $d_1 \geq 3$ and $e\geq 8$, this
implies $m_1 \geq 7$, but then $m_1d_1 \geq 21$, contradicting again $d(Z) = \Sigma m_id_i \leq 16$.
\par
Now assume $d_1\geq 5$. Since $m_1 \geq 2$ and $m_1d_1 \leq 16$, the only possibilities are: $m_1=2$,$5 \leq d_1 \leq 8$ and $d_i \leq 4$ for $i>1$, or $m_1=3,d_1=5$, furthermore in this case $Y_{red}=Y_1$ and $\tilde{Y}_1$ is subcanonical. In any case we can apply Proposition \ref{pmuc}, we just have to estimate $\delta$. Clearly for given $d_1$, $\delta \leq G(d_1)$, the maximal genus of a non-degenerated integral curve of degree $d_1$ in $\Pt$. Also if $d_1$ is even we can even take $\delta \leq G(d_1)-1$ because in this case curves of maximal genus are all complete intersections, hence Gorenstein. The corresponding values for $5\leq d_1\leq 8$ are $2,3,6,8$; in particular we may assume $\delta \leq d_1$. Arguing as above but using Prop. \ref{pmuc} instead of Prop. \ref{pmu}, with the bound $\delta \leq d_1$, yields the result. \epf

\subsubsection{$Z_{red}$ doesn't contain any irreducible component of degree three.}

Let's us first recall some well known facts about irreducible, non-degenerated, degree three
surfaces in $\Pq$. If $T$ is such a surface then $T$ is either a cubic scroll, $T'$, or a cone,
 $T''$, over a twisted cubic. 
\par
A cubic scroll, $T'$, is isomorphic to $\textbf{P}(\oc \oplus \oc (-1))$;
we have $Pic(T') = C_0{\textbf Z} \oplus f{\textbf Z}$ where $C_0^2 = -1, C_0.f = 1, f^2 =0$. 
The hyperplane system is $|C_0+2f|$, while  the curves of $\delta =|C_0+f|$ are conics. 
The linear system $\delta$ is $\infty ^2$ and the base locus of $\delta$ is empty. If $K \in \delta$
, we will denote by $\Pi _K$ the plane spanned by $K$. The planes $\Pi _K$ fill up $\Pq$ and two
such general planes intersect in one point.
\par
If $T''$ is a cone, then $T''$ is the image of ${\textbf P} (\oc \oplus \oc (-3))$ through $|C_0+3f|$,
where $C_0^2=-3, C_0.f=1, f^2 =0$. The only conics on $T''$ are the images of the curves in 
$|C_0+2f|$ which are pairs of rulings. If $K$ is a pair of rulings we will still denote by
$\Pi _K$ the plane they span. There are $\infty ^2$ such planes which fill up $\Pq$ and two such
general planes intersect in the vertex, $v$, which is the base locus of the linear system $\delta$
of the conics on $T''$.

\begin{proposition}
\label{compd=3}
With notations as in \ref{not2}, $Z_{red}$ doesn't contain any irreducible component of degree
three.
\end{proposition}

\textbf{Proof.}
Let $\tilde{T}$ denote an irreducible component of degree three of $Z_{red}$. By Proposition 
\ref{ZrinSing}, $\tilde{T} \subset Sing(\Sigma)$. Denote by $T=\tilde{T} \cap \Pq$ the section
with a general $4$-dimensional linear subspace. Also set $S=X \cap \Pq$ and $\Sc = \Sigma \cap
\Pq$. So we have the following situation:
\par
(i) $T \subset Sing(\Sc)$, the smooth surface $S$ is contained in $\Sc$.
\par
(ii) $dim(T \cap S)=0$ 
\par
The last assertion follows from the fact that $dim(X \cap Z)=2$ (see Lemma \ref{l2}).
\begin{enumerate}
\item First assume that if $T$ is a cone then $S$ doesn't pass through the vertex $v$ of the cone.
\par
Since $dim(S \cap T)=0$ and since the base locus of $\delta$ is empty if $T$ is a scroll (resp.
$=\{v\}$ if $T$ is a cone), if $K$ is sufficiently general in $\delta$ then $K \cap S= \emptyset$.
It follows that $dim(S \cap \Pi _K)=0$. Of course the general plane $\Pi _K$ is not contained in
$\Sc$, so $\Sc \cap \Pi _K$ is a plane quintic curve. This quintic curve contains twice the conic
$K$ (because $T \subset Sing(\Sc)$), thus $\Sc \cap \Pi _K=2K \cup L$. Since $(S \cap \Pi _K)
\subset (\Sc \cap \Pi _K)$ and since $S \cap K=\emptyset$, it follows that $(S\cap \Pi _K)$ is 
contained in the line $L$. This implies that $S$ is degenerated, which is absurd.
\item Assume now that $T$ is a cone and that $S$ passes through the vertex of the cone.
\par
We go back to $\Psx$. We are in the following situation: any we cut with a general $\Pq$, 
$\tilde{T} \cap \Pq$ is a cone of vertex $v$ and $v \in X \cap \Pq$. It follows that
$\tilde{T}$ is a cone of vertex a plane $\Pi$ over a twisted cubic and that $\Pi \subset X$.
Since $\Pi \subset X \cap Z_{red} \subset X \cap Z$. By Lemma \ref{pidsXZ} below we get 
$d \leq 4e+21$.
This implies (see \ref{not2}) $d(Z) \leq -e+11$. Since $e \geq 8$ (\cite{HS}), we get
$d(Z) \leq 3$ so $Z = \tilde{T}$ but this contradicts the fact that $Z$ is non-reduced (see Lemma
\ref{mi>1}). \epf
\end{enumerate}

\begin{lemma}
\label{pidsXZ}
Let $X \subset \Psx$ be a smooth codimension two subvariety which is not a complete intersection. 
Let $E$ be the rank two vector bundle 
associated to $X$: 
$$
0 \to \oc \to E \to \ic _X (e+7) \to 0
$$
Assume $s <e+7$, then a section of $\ic _X(s)$ yields a section of $E(-e-7+s)$ vanishing in 
codimension two:
$$
0 \to \oc \stackrel{\sigma}{\to}  E(-e-7+s) \to \ic _Z(-e-7+2s) \to 0 
$$
If $Z\cap X$ contains a two-dimensional plane, then $d \leq (s-1)e+5s-4$.
\end{lemma}

\textbf{Proof.}
Assume $X \cap Z$ contains a two dimensional plane, $\Pi$.
\par \noindent
$\underline {Claim:}$ There exists $k \geq 0$ such that $E_{\Pi}(-e-8+s-k)$ has a section vanishing 
in codimension two.
\par \noindent
The restriction of $\sigma$ to $\Pi$ vanishes identically since $\Pi \subset Z=(\sigma)_0$. Let $H$ 
be a general $\Pt$ containing $\Pi$ such that $H$ is not contained in $Z$ (we can always find such 
an $H$ because there are $\infty ^3$ $\Pt$'s containing $\Pi$). The restriction $\sigma |H$ vanishes 
on the divisor $\Pi$, dividing out by the equation of $\Pi$, we get $h^0(E_H(-e-8+s))\neq 0$. 
Repeating if necessary this argument, we reach the conclusion of the claim.
\par \noindent
Since $\Pi \subset X$ we have the exact sequence:
$$
0 \to N_{\Pi , X} \to N_{\Pi} \to N_X |\Pi \to 0
$$
Observe that $N_X = E_X$ the restriction of $E$ to $X$. So the exact sequence of normal bundles 
reads like:
$$
0 \to N_{\Pi , X} \to 4.\oc _{\Pi}(1) \stackrel{f}{\to} E_\Pi \to 0
$$
Twisting by $\oc (-1)$ we get that $E_{\Pi}(-1)$ is generated by $4$ global sections; moreover a 
general section in $Im(H^0(f))$ has a smooth zero-locus of codimension two, $\Gamma$.
$$
0 \to \oc _{\Pi} \to E_{\Pi} (-1) \to \ic _{\Gamma} (e+5) \to 0
$$
Since $E_\Pi (-1)$ is generated by $4$ global sections, we get that $\ic _{\Gamma} (e+5)$ is 
generated by $3$ global sections:
$$
0 \to N_{\Pi ,X}(-1) \to 3.\oc _{\Pi} \to \ic _{\Gamma} (e+5) \to 0
$$
It follows that $\Gamma$ is contained in an irreducible (actually smooth) curve of degree $e+5$. 
On the other hand, from the claim we have $h^0(E_{\Pi} (-e-8+s)) \neq 0$; this implies, 
since $-e-7+s <0$ by our assumption, that $h^0(\ic _{\Gamma} (s-2)) \neq 0$. Hence $\Gamma$ is 
contained in a complete intersection of type $(s-2, e+5)$, and therefore $d(\Gamma) \leq (s-2)(e+5)$.
We have $d(\Gamma) = c_2(E_{\Pi} (-1)) = -c_1(E)+c_2(E)+1 =-e-7+d+1 = d-e-6$. In conclusion we have:
$d \leq (s-2)(e+5)+e+6$. \epf

\subsubsection{$Z_{red}$ doesn't contain any irreducible component of degree at least four.}

So far we have seen that every irreducible component of $Z_{red}$ has degree at least four
and that $Z_{red} \subset Sing(\Sigma)$.
\par
Let $\Pt$ be a general three-dimensional linear subspace of $\Psx$ and set: $\tilde{Y} = Z_{red}\cap
 \Pt$,
$F=\Sigma \cap \Pt$, $C=X\cap \Pt$. We have $\tilde{Y} \subset Sing(F)$, $C\subset F$ and 
$C\cap Sing(F)=
\emptyset$ (the last assertion follows from the fact that $dim(X\cap Sing(\Sigma))=2$, see Lemma
\ref{l2}.
\par
Let $f:\tilde{F} \to F$ be a desingularization of $F$ and set $\oc _{\tilde{F}}(H))=f^*(\oc _F(1))$.

\begin{lemma}
\label{Flin}
With notations as above, $h^0(\oc _{\tilde{F}}(H))=4$ (i.e. $F$ is ''linearly normal").
\end{lemma}

\textbf{Proof.}
If $h^0(\oc _{\tilde F}(H)) >4$ then $f$ factors through $\Pq$. Since $C$ is smooth and since 
$C \cap Sing(F)=\emptyset$, $C$ would be the isomorphic image of a curve of degree $d$ in $\Pq$, 
but this contradicts the linear normality of $C$ (see Lemma \ref{Clin}). \epf 

The following lemma will conclude the proof:

\begin{lemma}
\label{Fnlin}
Let $F \subset \Pt$ be an irreducible quintic surface. Assume $Sing(F)$ contains a (reduced)
curve $\tilde{Y}$ such that every irreducible component of $\tilde{Y}$ has degree at least four. 
Then $F$ is
rational or ruled over an elliptic curve. In any case $F$ is not ''linearly normal".
\end{lemma}

\textbf{Proof.}
Let $Y$ be an irreducible component of $\tilde{Y}$; $Y$ is an integral curve of degree at least
four. Also $deg(Y) \leq 6$ since $Y$ is contained in the singular locus of an irreducible
quintic surface. We may also assume that $Y$ is not planar.
\begin{enumerate}
\item Assume $deg(Y)=4$.
\par
If $p_a(Y)=0$ then $Y$ is a smooth rational quartic and lies on a quadric $Q$ which is the surface 
of trisecants to $Y$,since $Y \subset Sing(F)$, we should have $Q \subset F$ which is absurd.
\par
So we may assume $p_a(Y)=1$. In this case $Y$ is the complete intersection of two quadrics. 
If $p \in F\setminus Y$ then there exists a quadric, $Q_p$ of the pencil 
${\textbf P}(H^0(\ic _Y(2)))$ passing through $p$. The complete intersection $F \cap Q_p$ links 
$2Y$ to a conic $K_p$. 
\par \noindent
If $K_p$ is irreducible for general $p$ then $F$ is rational (look at the morphism 
$F \setminus Y \to {\textbf P}^1 \simeq {\textbf P}(H^0(\ic _Y(2)))$ whose general fiber is $K_p$).
\par \noindent
If $K_p$ is reducible for every $p$, then $F$ is ruled and through any point of $Y$ there pass 
two rulings (which form a conic $K_p$) of $F$. Let $H$ be the plane spanned by such a degenerated
conic. Then $H \cap F$ must contain a plane curve which is not a ruling, this curve has
degree at most three, hence it is of geometric genus at most one. This implies that $F$ is 
rational or ruled over an elliptic curve.
\item Assume $deg(Y)=5$.
\par
If $p_a(Y)=2$ then $Y$ is contained in a quadric surface, $Q$, which is the surface of trisecants to
$Y$; we should have $Q \subset F$, which is absurd.
\par
If $p_a(Y)=1$, then $Y$ has at most one singular point and through a general point $p\in Y$ there
is a trisecant to $Y$ (not passing through the singular point, if any) through $p$. Since $Y \subset
Sing(F)$ these trisecants are contained in $F$. It follows that the general point of $F$ is
contained in a unique such trisecant to $Y$. This shows that $F$ is ruled. Since the plane section
of $F$ has geometric genus $\leq 1$, $F$ is rational or ruled over an elliptic curve. The same
argument applies if $p_a(Y)=0$.
\item Finally assume $deg(Y)=6$.
In this case the plane section of $F$ is a rational curve and $F$ is rational.
\end{enumerate}
\par
To conclude it remains to show that if $F$ is rational or ruled over an elliptic curve, then
$F$ is not "linearly normal".
\par
Consider the exact sequence
$$
0 \to \oc _{\tilde F} \to \oc _{\tilde F}(H) \to \oc _H(H) \to 0
$$
since $h^1(\oc _{\tilde F})\leq 1$, $H^2 = 5$ and since $H$ is a smooth curve of genus $1$ if $F$
is ruled over an elliptic curve (resp. of genus $\leq 2$ if $F$ is rational),
we get: $h^0(\oc _{\tilde F}(H))\geq 5$. \epf
\bigskip

\begin{corollary}
\label{compd>3}
With notations as in \ref{not2}, $Z_{red}$ has an irreducible component of degree at most three.
\end{corollary}

\textbf{Proof.}
Follows directly from Lemma \ref{Flin} and Lemma \ref{Fnlin}. \epf 

\subsubsection{Conclusion.}

We can now state the main result of this paper:

\begin{theorem}
\label{thm}
Let $X \subset \Psx$ be a smooth subvariety of codimension two. If $h^0(\ic _X(5)) \neq 0$, then
$X$ is a complete intersection.
\end{theorem}

\textbf{Proof.}
If $h^0(\ic _X(4)) \neq 0$, this follows from Ran's theorem (\cite{Ran}). If $X$ lies on an
irreducible quintic hypersurface, the result follows from Prop. \ref{compd=1}, \ref{compd=2}, 
\ref{compd=3} and Corollary \ref{compd>3}. \epf

\section{Subcanonical threefolds in $\Pcq$ lying on a quartic hypersurface.}

In this section we show how the previous methods apply to prove:

\begin{theorem}
\label{P5}
Let $X \subset \Pcq$ be a smooth, subcanonical threefold. If $s\leq 4$ then $X$ is a complete intersection.
\end{theorem}

\textbf{Proof.} First of all we can assume $e\geq 7$: indeed, by \cite{Ball_Ch}, Prop. 9, 10, if $e\leq 2$, then $X$ is a complete intersection; if $e\geq 3$ the rank two vector bundle associated to $X$ satisfies $c_1^2\geq 4c_2$ and we apply \cite{Ho}, Cor.7.3. 
\par
As in \ref{not2} we see that $X$ is bilinked on the quartic hypersurface $\Sigma$ to a subscheme $Z$ with $d(Z)=d-4e-8$ and $\omega _Z \simeq \oc _Z(-e-4)$. As in Lemma \ref{mi>1}, every irreducible component of $Z_{red}$ appears with multiplicity in $Z$. By Lemma \ref{lsd}, $d(Z)\leq 9$.
\par
Assume first that $Z_{red}$ contains an irreducible component, $L$, of degree one. Arguing as in  the proof of Prop. \ref{compd=1}, we get $S \subset \Pq$ with $dim(Sing(S))\leq 0$ and $S \subset L\cup F$ where $F$ is an hypersurface of degree three. Since $S$ has at most isolated singularities and is locally Cohen-Macaulay, by \cite{Ha2}, $S$ is irreducible, hence we may assume $S \subset F$. We may assume $F$ irreducible (otherwise we conclude considering a general hyperplane section of $S$). Since $q(X)=0$ for a smooth volume in $\Pcq$, a curve which is a space section of $X$ is linearly normal (cf Lemma \ref{Clin}). It follows that $F$ is normal and we conclude with \cite{Koe}. From now on we may assume that every irreducible component of $Z_{red}$ has degree at least two. On the other hand every such irreducible component has degree at most four, hence its general space section is a Gorenstein curve. Applying Prop. \ref{pmu}, we deduce that $Z_{red} \subset Sing(\Sigma)$ (cf Prop. \ref{ZrinSing}). Observe that an irreducible quartic surface in $\Pt$, containing an irreducible curve of degree at least two in its singular locus, is rational, hence one has that such a surface is not "linearly normal" (in the sense of Lemma \ref{Flin}), this contradicts the linear normality of a general space section of $X$. \epf

\section{Low degrees in $\Psx$.}

Let's recall the following result of Holme and Schneider (see \cite{HS}, Cor. 6.3 and the proof 
of Cor. 6.2):
\bigskip

\begin{proposition} Let $X \subset \Psx$, be a smooth subvariety of codimension
two, of degree $d \leq 73$. Then $X$ is a complete intersection unless $X$ is numerically a
complete intersection of type $(a,b)$, $a \leq b$, for some $(a,b)$ in the table:
\bigskip

\begin{tabular}{|c||c|c|c|c|c|c|c|c|}
\hline
(a,b) & (7,8) & (6,10) & (7,9) & (8,8) & (6,11) & (7,10) & (8,9) & (6,12) \\ \hline
d & 56 & 60 & 63 & 64 & 66 & 70 & 72 & 72 \\ \hline
e & 8 & 9 & 9 & 9 & 10 & 10 & 10 & 10 \\ \hline
\end{tabular}
\end{proposition}
\bigskip

\begin{proposition}
\label{Pnci}
Let $X \subset \Psx$ be numerically a complete intersection $(a,b),a \leq b$.
\begin{enumerate}
\item If $a \leq 6$ then $X$ is a complete intersection.
\item If $b \leq a+2$ and $X$ is not a complete intersection, then $s \leq a-2$ and $a \geq 8$.
\item If $(a,b)=(7,10)$ or $(8,8)$ then $s \geq 7$.
\item If $a=7$ and $b \leq 10$, then $X$ is a complete intersection.
\end{enumerate}
\end{proposition}

\textbf{Proof.}
1) If $X$ is not a complete intersection then $s \leq a-1$ (see the proof of \ref{nci}), but by Theorem \ref{thm}, $6 \leq s$, so $a \geq 7$.
\par
2) Assume $s=a-1$ and $b \leq a+2$. Consider the rank two vector bundle associated to $X$:
$$
0 \to \oc \to E \to \ic _X(a+b) \to 0
$$
then $E(-b-1)$ has a section vanishing in codimension two:
$$
0 \to \oc \to E(-b-1) \to \ic _Z(a-b-2) \to 0
$$
We have $\omega _Z \simeq \oc _Z(a-b-9)$ and $d(Z)=b-a+1$. Since $b \leq a+2$, $d(Z) \leq 3$. Arguing as in the proof of \ref{mi>1}, we see that every irreducible component of $Z_{red}$ appears with a non-reduced structure in $Z$. It follows that $Z$ is a l.c.i. multiple structure on a linear space, of multiplicity $d(Z) \leq 3$. By \cite{Ma}, $Z$ is a complete intersection and $E$ splits, contradiction. So $s<a-1$. Since by Theorem \ref{thm} $s \geq 6$, we get $a \geq 8$.
\par
3) Now assume $(a,b)=(7,10)$ or $(8,8)$ and $s \leq 6$. By Thm \ref{thm} we may assume $s=6$. With notations as in 2), we have $d(Z)=4$ and $\omega _Z \simeq \oc _Z(-a-b+5)$. Again we see that every irreducible component of $Z_{red}$ appears in $Z$ with a non-reduced structure. If $Z_{red}$ contains an irreducible component of degree one then $Z$ is either: a) a multiplicity four l.c.i. structure on a linear space, or b) a double structure on $L_1 \cup L_2$ where $L_i$ is a linear space. In case a), by \cite{Ma}, $Z$ is a complete intersection and we are done. In case b) we claim that $L_1 \cup L_2$ is contained in a hyperplane. Indeed, otherwise cutting with a general $\Pq$, $Z\cap \Pq$ would have support on two planes meeting in one point, by \cite{Ha2}, $Z\cap \Pq$ cannot be loc.C.M., in contradiction with the fact that $Z\cap \Pq$ is l.c.i. So in any case we may assume that $Z_{red}$ is a quadric, $Q$, contained in an hyperplane. 
\par
Let $\Sigma$ denote the sextic hypersurface containing $X$. Then $X$ and $Z$ are bilinked on $\Sigma$ and we claim that: (i) $dim(X \cap Sing(\Sigma))=2$, (ii) $Z_{red} \subset Sing(\Sigma)$.
\par
If (i) is not true by Lemma \ref{lsd}, since $h^1(\ic _X(1))=0$, we get $d \leq 4(e+5)+1$, which is impossible in our cases. For (ii) we argue as in Prop. \ref{ZrinSing}: consider $Y=Z\cap \Pt$ and assume $Y=Y_1 \cup Y_2$ where $Y_i$ is a double structure on the line $R_i$. Since $\omega _Y \simeq \oc _Y(-a-b+8)$, we get $2p_a(Y_i)-2 \leq 2(-a-b+8)$. If $R_i$ is not contained in $Sing(\Sigma \cap \Pt)$, then $2(-a-b+8) \geq 2\mu (1,0,6,2) -2$ and one checks that in our cases this inequality is not satisfied. The same argument applies if $Y$ is a double structure on a smooth conic.
\par
Now consider the hyperplane $H \subset \Psx$ containing the quadric $Q=Z_{red}$. We have $V=(X\cap H) \subset (\Sigma \cap H)=2Q \cup F$. Since
$dim(X \cap Z)=2$, it follows that $V \subset F$, i.e. $h^0(\ic _V(2)) \neq 0$. From the exact sequence $0 \to \ic _X(1) \to \ic _X(2) \to \ic _V(2) \to 0$, since $h^1(\ic _X(1))=0$, we get $h^0(\ic _X(2))\neq 0$, a contradiction.
\par
4) Follows from 2) and 3). \epf

\begin{theorem}
\label{thm2}
Let $X \subset \Psx$ be a smooth subvariety of codimension two and degree $d$. If $d \leq 73$,
then $X$ is a complete intersection.
\end{theorem}

\textbf{Proof.}
Look at the table above. If $a \leq 6$, we apply A) of Prop. \ref{Pnci}. If $a=7$, and $b\leq 10$ we conclude with C) of Prop. \ref{Pnci}. Finally if $(a,b)=(8,8)$, by C) of Prop. \ref{Pnci}, $s=7$. In this case we have $d(Z)=1$ and we easily conclude.
\par
It remains to do the case $(8,9)$. We may assume $6\leq s \leq 7$. If $s=7$ then $E(-10)$ has a section vanishing along a codimension two subscheme, $Z$, of degree $c_2(E(-10))=2$; it follows that $Z$ is a complete intersection.
\par
Now assume $s=6$, in this case $E(-11)$ has a section vanishing along a codimension two subscheme, $Z$, of degree $6$ with $\omega _Z \simeq \oc _Z(-12)$.
\par
The subscheme $Z$ is bilinked to $X$ on the sextic hypersurface, $\Sigma$, containing $X$. Using Lemma \ref{lsd} we see that $dim(X\cap Z)=dim(X\cap Sing(\Sigma))=2$ (cf Lemma \ref{l2}). We claim that at least one irreducible component of $Z_{red}$ is contained in $Sing(\Sigma)$.
\par \noindent
Indeed with notations as in \ref{not2}, consider $Y=Z\cap \Pt$, $Y=\tilde{Y}_1 \cup ... \tilde{Y}_r$, where $\tilde{Y}_i$ is a multiplicity $m_i$ structure on the integral curve $Y_i$ of degree $d_i$. We have $m_i \geq 2$ (see Lemma \ref{mi>1}) and $\Sigma m_id_i=6$. Of course $d_i\leq 3$ so that $Y_i$ is Gorenstein and we can apply Prop. \ref{pmu}. Arguing as in the proof of Prop. \ref{ZrinSing}, we have $2p_1-2 \leq -9m_1d_1$ ($p_1=p_a(\tilde{Y}_1$)), but if $Y_1$ is not contained in $Sing(F)$ ($F=\Sigma \cap \Pt$), we must also have $-9m_1d_1 \geq 2\mu (d_1,g_1,6,m_1)-2$. It is easily seen that this inequality can be satisfied only if $d_1=1,g_1=0$ and $m_1 \geq 3$. If $m_1\geq 5$ then $Z$ is a l.c.i. multiplicity $6$ structure on a linear subspace, by \cite{Ma}, $Z$ is a complete intersection and we are done. If $m_1=4$ then $m_2=2,d_2=1$ and $Y_2 \subset Sing(F)$. It remains the case $m_1=m_2=3,d_1=d_2=1$. In this case we can see $Y$ as a triple structure on the (reducible) conic $Y_1\cup Y_2$; applying again Prop. \ref{pmu} (which is possible since $Y_1\cup Y_2$ is Gorenstein), we get $Y_1\cup Y_2 \subset Sing(F)$.
\par 
In conclusion $Z_{red}$ has an irreducible component of degree $d',1\leq d' \leq 3$, contained in $Sing(\Sigma)$.
\par
If $d'=1$, arguing like in the proof of Prop. \ref{compd=1} , we get a smooth surface $S \subset \Pq$ with $h^0(\ic _S(4))\neq 0$. It follows (Prop. \ref{lsd}), that $72=d \leq 4(5+10)+1$, which is absurd.
\par
If $d'=2$, we get an hyperplane section $V=X\cap H$, with $h^0(\ic _V(2))\neq 0$ (cf proof of Prop. \ref{compd=2}). Since $h^1(\ic _X(1))=0$ by Zak's theorem, this implies $h^0(\ic _X(2))\neq 0$, a contradiction.
\par
Finally we are left with the case $d'=3$ where $Z$ is a double structure on an integral subvariety of degree three. We argue like in Prop. \ref{compd=3} so let's denote by $T$ the intersection of $Z_{red}$ with a general $\Pq$. Like in the first part of the proof of Prop. \ref{compd=3}, we see that in the case $T$ smooth (or $S$ not passing through the vertex of $T$ if $T$ is a cone), $S$ has a plane section, $S\cap \Pi _K$, contained in a conic. Let $C=S\cap H$ where $H$ is a general hyperplane through $\Pi _K$. By Lemma \ref{Clin}, $h^1(\ic _C(1))=0$, it follows that $h^0(\ic _C(2))\neq 0$, this in turn implies $h^0(\ic _S(2))\neq 0$ and from this it follows that $S$ is a complete intersection.
\par
We are left with case $2$ of the proof of Prop. \ref{compd=3} i.e. $Z_{red}$ is a cone of vertex a plane $\Pi$ over a twisted cubic and $\Pi \subset X$. Unfortunately Lemma \ref{pidsXZ} is insufficient to conclude, so we try to improve the argument. First observe that a general $\Pq$ through $\Pi$ intersects $Z_{red}$ in $\tilde{\Pi}$, the first infinitesimal neighbourhood of $\Pi$ in $\Pq$. Thus $Z$ is a double structure on $\tilde{\Pi}$. As in the proof of Lemma \ref{pidsXZ}, consider a general $\Pt$ through $\Pi$, the intersection $\Pt \cap Z$ contains the divisor $2\Pi$; it follows that $h^0(E_{\Pi}(-e-9+s))\neq 0$ and going on with the argument we get $d\leq (s-3)(e+5)+e+6$, which in our case is impossible.  \epf

Author's address: Dipartimento di Matematica, via Machiavelli, 35. 
\par
44100 Ferrara (Italy)
\par
emails: phe@dns.unife.it (Ph.E.), d.franco@dns.unife.it (D.F.)

\end{document}